\documentclass[a4paper]{IEEEtran}

\usepackage{balance}
\usepackage{color}

\IEEEoverridecommandlockouts
\usepackage{graphicx,epsfig}

\usepackage{algorithmic}
\usepackage{algorithm}

\usepackage[cmex10]{amsmath}
\usepackage{url}
\usepackage{multirow}

\newtheorem{theorem}{Theorem}[section]

\title{An error estimate of Gaussian Recursive Filter\\ in 3Dvar problem}
\author{
Salvatore Cuomo$^{(1)}$, Raffaele Farina$^{(2)}$, Ardelio Galletti$^{(3)}$, Livia Marcellino$^{(3)}$\\[0.3cm]
$^{(1)}$University of Naples Federico II\\
Department of Mathematics and Applications ''R. Caccioppoli'', Italy\\
Email: salvatore.cuomo@unina.itl\\[0.3cm]
$^{(2)}$Centro Euro-Mediterraneo sui Cambiamenti Climatici\\
CMCC, Italy\\
Email: raffaele.farina@cmcc.it\\[0.3cm]
$^{(3)}$University of Naples "Parthenope"\\
Department of Science and Technology, Italy\\
Email: \{ardelio.galletti, livia.marcellino\}@uniparthenope.it\\
}

\begin{document}
\maketitle              

\begin{abstract}
Computational kernel of the three-dimensional variational data
assimilation (3D-Var) problem is a linear system, generally solved
by means of an iterative method. The most costly part of each
iterative step is a matrix-vector product with a very large
covariance matrix having Gaussian correlation structure. This
operation may be interpreted as a Gaussian convolution, that is a very
expensive numerical kernel. Recursive Filters (RFs) are a
well known way to approximate the Gaussian convolution and are  intensively applied
in the meteorology, in the oceanography and in forecast models.
In this paper, we deal with an oceanographic 3D-Var
data assimilation scheme, named OceanVar, where the linear system is solved by using the Conjugate Gradient (GC) method by replacing, at each
step, the Gaussian convolution with RFs. Here we give
theoretical issues on the discrete convolution approximation
with a first order (1st-RF) and a third order  (3rd-RF) recursive filters.
Numerical experiments confirm given error bounds and show the benefits,
in terms of accuracy and performance, of the 3-rd RF.
\end{abstract}

\section{Introduction}
In recent years, Gaussian filters have assumed a central role in image filtering and techniques for accurate measurement \cite{Witkin}. The implementation of the Gaussian filter in one or more dimensions has typically been done as a convolution with a Gaussian kernel,  that leads to a high computational cost in its practical application. Computational efforts to reduce the Gaussian convolution complexity  are discussed in \cite{Haglund, Vliet3}. More advantages may be gained by employing a \emph{spatially recursive filter}, carefully constructed to mimic the Gaussian convolution operator.\\
Recursive filters (RFs) are an efficient way of achieving a long impulse response, without having to perform a long convolution.  Initially developed in the context of time series analysis \cite{Dahlquist},  they are extensively used as computational kernels for numerical weather analysis, forecasts \cite{ Lorenc, Purser, Weaver}, digital image processing \cite{Deriche, Vliet}. Recursive filters with higher order accuracy are very able to accurately  approximate a Gaussian
convolution, but they require more operations.\\
In this paper, we investigate how the RF mimics  the Gaussian convolution
in the context of variational data assimilation analysis.
Variational data assimilation (Var-DA) is popularly used to combine
observations with a model forecast in order to produce a \emph{best}
estimate of the current state of a system and enable accurate
prediction of future states. Here we deal with the three-dimensional
data assimilation scheme (3D-Var), where the estimate minimizes a
weighted nonlinear least-squares measure of the error between the
model forecast and the available observations.
The numerical problem is to minimize a cost function by means of an
iterative optimization algorithm. The most costly part of each
 step is the multiplication of some
grid-space vector by a covariance matrix that defines the error on
the forecast model and observations. More precisely, in 3D-Var
problem this operation may be interpreted as the convolution of a
covariance function of background error with the given forcing
terms.\\
Here we deal with  numerical aspects of an oceanographic
3D-Var scheme, in the real scenario of OceanVar. Ocean data assimilation is a
crucial task in operational oceanography and the computational
kernel of OceanVar software is a linear system resolution by means
of the Conjugate Gradient (GC)  method, where the iteration matrix
is relate to an errors covariance matrix, having a
Gaussian correlation structure.\\
In \cite{Dobricic}, it is shown that a computational advantage
can be gained by employing a first order RF that mimics the required Gaussian convolution.
Instead, we use the 3rd-RF to compute numerically the Gaussian convolution,
as how far is only used in signal processing \cite{Young}, but only recently used in the field of Var-DA problems.\\
In this paper we highlight the main sources of error,
introduced by these new numerical operators.
We also investigate the real benefits, obtained by using 1-st and 3rd-RFs,
through a careful error analysis. Theoretical aspects are confirmed by
 some numerical experiments. Finally, we report
results in the case study of the OceanVar software.\\
The rest of the paper is organized as follows. In the next section we
recall the three-dimensional variational data assimilation problem
and we remark some properties on the conditioning for this problem.
Besides, we describe our case study: the OceanVar problem and its
numerical solution with CG method. In section III, we
introduce the $n$-th order recursive filter and how it can be
applied to approximate the discrete Gaussian convolution. In section IV, we estimate
the effective error, introduced at each iteration of the CG method, by using 1st-RF and 3rd-RF
instead of the Gaussian convolution. In section V, we
report some experiments to confirm our theoretical study, while the
section VI concludes the paper.

\section{Mathematical Background}
The aim of a generic variational problem ({\sc Var} problem) is to
find a best estimate $x$, given a previous estimate $x_b$ and a
measured value $y$. With these notations, the {\sc Var} problem is
based on the following regularized constrained least-squared
problem:
\[
\min_{x } J(x)
\]
where $x$ is defined in a grid domain $D$. The objective function
$J(x)$ is defined as follows:
\begin{equation}\label{LS pb}
J(x) = \| y- \mathcal H(x)\|^{2} + \lambda R(x, x_b)
\end{equation}
where measured data are compared with the  solution obtained from a nonlinear model given by $\mathcal H(x)$.\\
In (\ref{LS pb}), we can recognize a quadratic data-fidelity term, the first term and the general regularization term (or penalty term), the second one. When \(\lambda=1\) and the regularization term can be write as:
\[
 R(x, x_b)= \| x-x_b\|^{2}
\]
we deal with a three-dimensional variational data assimilation
problem (3D-Var DA problem). The purpose is  to
find an optimal estimate for a vector of states \(x_t\) (called the
analysis) of a generic system \(S\), at each time \(t \in T=
\{0,..,n\}\) given:
\begin{itemize}
\item a prior estimate vector \(x_t^b\) (called the background) achieved by numerical solution of a forecasting model \( \mathcal L_{t-1,t}(x_{t-1})=x_t^b\), with error \( \delta x_t = x_t^b - x_t\);
\item a vector \(y_t\) of observations, related to the nonlinear model by \( \delta y_t \) that is an effective measurement error:
\[
y_t= H(x_t)+\delta y_t.
\]
\end{itemize}
At each time t, the errors \( \delta x_t \) in the background and the errors \( \delta y_t\) in the observations are assumed to be random with mean zero and covariance matrices \({\bf B}\) and \({\bf R}\), respectively. More precisely, the covariance \({\bf R}=<\delta y_t , \delta y_t^T>\) of observational error is assumed to be diagonal, (observational errors statistically independent). The covariance \({\bf B}=<\delta x_t , \delta x_t^T>\) of  background error is never assumed to be diagonal
as justified in the follow. To minimize, with respect to \(x_t\) and for each \(t \in T\), the
problem becomes:
\begin{equation}\label{DA pb}
\min_{x_t \in D}J(x_t) = \min_{x_t \in D} \{  \frac{1}{2}     \|
y_t-H(x_t)\|_{\bf R}^{2} + \frac{1}{2}\| x_t-x_t^b\|_{\bf B}^{2} \}
\end{equation}
In explicit form, the functional cost of (\ref{DA pb}) problem can
be written as:
\begin{equation}\label{DA1 pb}
 \begin{array}{l}
J(x_t) = \frac{1}{2}  (y_t-H(x_t))^T {\bf R}^{-1} (y_t-H(x_t)) +\\
\\
\hspace{1cm} + \frac{1}{2} (x_t-x_t^b)^T {\bf B}^{-1}(x_t-x_t^b)
\end{array}
\end{equation}
It is often numerically convenient to approximate the effects on
\(H(x_t)\) of small increments of \(x_t\), using the linearization
of \(H\). For small increments \( \delta x_t\), follows
\cite{Lorenc2}, it is:
\[
H(x_t) \simeq H(x_t^b) + {\bf H} \delta x_t
\]
where the linear operator \({\bf H}\) is the matrix obtained by the
first order approximation of the Jacobian of \(H\) evaluated at \(x_t^b\).\\
Now let  \(d_t= y_t-H(x_t^b)\) be the \emph{misfit}. Then the
function \(J\) in  (\ref{DA1 pb})   takes the following form in the
increment space:
\begin{equation}\label{DA2 pb}
\!\!\!\begin{array}{l}\!J(\delta x_t) \!=\! \frac{1}{2} (d_t \!-\!
{\bf H} \delta x_t)^T {\bf R}^{-1}\!(d_t \!-\! {\bf H} \delta x_t)\!
+\!\frac{1}{2} \delta x_t^T {\bf B}^{-1} \delta x_t\!\!\!\!
\end{array}
\end{equation}
At this point, at each time \(t\), the minimum of (\ref{DA2 pb}) is
obtained by requiring  $\nabla
J = 0$. This gives rise to the linear system:
\[
({\bf B}^{-1}+{\bf H}^T{\bf R}^{-1}{\bf H}) \delta x_t = {\bf H} ^T
{\bf R}^{-1} d_t
\]
or equivalently:
\begin{equation}\label{system1}
 (I+{\bf B}{\bf H}^T{\bf R}^{-1}{\bf H}) \delta x_t = {\bf B}{\bf H} ^T {\bf R}^{-1} d_t
\end{equation}
For each time $t=0,...,n$, iterative methods, able to converge
toward a practical solution, are needed to solve the linear system
(\ref{system1}). However this problem, so as formulated, is
generally very ill conditioned.  More precisely, by following
\cite{Haben}, and assuming that
\begin{equation}
{\bf \Psi} = {\bf H}^T{\bf R}^{-1}{\bf H}
\end{equation} \label{Lambda_M}
is a diagonal matrix, it can be proved that the conditioning of $I+{\bf B}{\bf \Psi}$ is strictly related to the conditioning of the matrix \({\bf B}\) (the covariance matrix). In general, the matrix \({\bf B}\) is a block-diagonal matrix, where each block is related to a single state of vector \(x_t\) and it is ill conditioned.\\
This assertion is exposed in  \cite{Nichols} starting from the
expression of  \({\bf B}\) for one-state vectors as:
\[ {\bf B}= \sigma_b^2 {\bf C} \]
 where \(\sigma_b^2\) is the background error variance and \({\bf C}\) is a matrix that denotes
 the correlation structure of the background error. Assuming that the correlation structure of matrix \({\bf C}\)
 is homogeneous and depends only on the distance
 between states and not on positions, an expression of \({\bf C}\) as a symmetric matrix
 with a circulant form is given; i. e. as a Toeplitz matrix. By means of a spectral analysis of its eigenvalues,
 the ill-conditioning of the matrix \({\bf C}\) is checked. As in \cite{Damore},  it follows that \({\bf B}\)  is ill-conditioned
 and the matrix \(I+{\bf B}{\bf \Psi}\), of the linear system
 (\ref{system1}),
 too. A well-known technique for improving the convergence of  iterative
methods for solving linear systems is to \emph{preconditioning} the system and thus reduce
the condition number of the problem.\\ In order to precondition the
system in (\ref{system1}), it is assumed that \({\bf B}\) can be
written in the form \({\bf B}={\bf V}{\bf V}^T\),  where  \({\bf
V}={\bf B}^{1/2}\) is the square root of the background error
covariance matrix \({\bf B}\). Because \({\bf B}\) is symmetric Gaussian,
\({\bf V}\) is uniquely defined as the symmetric (${\bf V^T}={\bf
V}$) Gaussian matrix such that ${\bf V^2}={\bf B}$. \\
As explained in \cite{Lorenc2}, the cost function (\ref{DA2 pb})
becomes: $$\!\!\!\!\begin{array}{l} J(\delta x_t) \!= \! \frac{1}{2}
(d_t \!- \!{\bf H} \delta x_t)^T {\bf R}^{-1}  (d_t \!- \! {\bf H}
\delta x_t) \!+\!\frac{1}{2}
\delta x_t^T ({\bf V}{\bf V^T})^{-1} \delta x_t \\[0.3cm]
\ \  \quad \!=\!  \frac{1}{2}  (d_t \!- \! {\bf H} \delta x_t)^T
{\bf R}^{-1}  (d_t \!- \! {\bf H} \delta x_t) + \frac{1}{2} \delta
x_t^T ({\bf V^T})^{-1} {\bf V}^{-1}\delta x_t
\end{array}
$$
Now, by using a new control variable \(v_t\), defined as  \(v_t=
{\bf V}^{-1} \delta x_t\), at each time \(t \in T\) and observing
that \(\delta x_t = {\bf V}v_t\) we obtain a new cost function:
 \begin{equation}\label{DA3 pb}
\widetilde {J}(v_t) = \frac{1}{2}  (d_t - {\bf H}{\bf V} v_t)^T {\bf
R}^{-1}  (d_t - {\bf H}{\bf V} v_t) + \frac{1}{2} v_t^T v_t.
\end{equation}
Equation (\ref{DA3 pb}) is said the \emph{dual problem} of equation
(\ref{DA2 pb}).  Finally, to minimize the cost
function $\widetilde {J}(v_t)$ in (\ref{DA3 pb}) leads to the new
linear system:
\begin{equation}\label{system2}
 (I+{\bf V}{\bf \Psi}{\bf V}) v_t = {\bf V}{\bf H}^T{\bf R}^{-1}d_t
\end{equation}
Upper and lower bounds on the condition number of the matrix $I+{\bf
V}{\bf \Psi}{\bf V}$  are shown  in \cite{Nichols}. In particular
it holds that:
\[
\mu(I+{\bf V}{\bf \Psi}{\bf V}) << \mu(I+{\bf B}{\bf \Psi}).
\]
Moreover, under some special assumptions, it can be proved that
$I+{\bf V}{\bf \Psi}{\bf V}$ is very well-conditioned ($\mu(I+{\bf
V}{\bf \Psi}{\bf V})<4$).

\subsection*{The OceanVar model}
\indent As described in \cite{Dobricic}, at each
time \(t \in T\), OceanVar software implements an oceanographic
three-dimensional variational DA scheme (3D Var-DA) to produce
forecasts of ocean currents for the
Mediterranean Sea. The computational kernel is based on the
resolution of the linear system defined in (\ref{system2}).
To solve it, the  Conjugate Gradient
(CG) method is used and a basic outline  is
described in  \textbf{Algorithm 1}.
\begin{algorithm}[h!]\label{CG}
    \caption{CG Algorithm}%
        \begin{algorithmic}[1]
            \STATE $k=0$; $\quad \mathbf x_0$, the initial guess;\\[0.5 mm]
            \STATE $\mathbf r_0 = \mathbf b-\mathbf A \mathbf x_0$;\\[0.5 mm]
            \STATE $\bf \rho_0 = \bf r_0$;\\[0.5 mm]
            \WHILE{$\big({ \|\mathbf r_{k}\|}/{\|\mathbf b\|} > \epsilon \ .and.\ k \leq n \big)$}
                \STATE $\mathbf q_k=\mathbf A \mathbf \bf \rho_k$;
                \STATE $ \alpha_k ={ (\mathbf{r}_k^T,\mathbf r_k)}/({\mathbf \rho_k, \mathbf q_k})$; $\quad \mathbf x_{k+1} = \mathbf x_{k} +\alpha_k \bf \rho_k$;\\[0.5 mm]
            \STATE ${\mathbf r_{k+1}} = {\mathbf r_k}-\alpha_k {\mathbf q_k}$; $\quad \beta_k ={({ \mathbf r_{k+1}^T}, \mathbf r_{k+1})}/{({ \mathbf r}_{k}^T,\mathbf r_{k})}$;\\[0.5 mm]
                \STATE $\mathbf \rho_{k+1}={\mathbf r}_{k+1}+\beta_{k} \mathbf \rho_{k}$; $\quad k=k+1$; \\[0.5 mm]
            \ENDWHILE
        \end{algorithmic}
\label{CGA}
\end{algorithm}

\noindent We focus our attention on step 5.: at each iterative step,
a matrix-vector product ${\bf A\, \rho}_k$ is required, where $$\bf
A=I+{\bf V}{\bf \Psi}{\bf V},$$  $\bf {\rho}_k$ is the residual at
step $k$ and ${\bf \Psi}$ depends  on the number of observations and
is characterized by a bounded norm (see \cite{Haben} for details).
More precisely, we look to the matrix-vector product
\[
\mathbf q_k=(I+{\bf V}{\bf \Psi}{\bf V})\bf {\rho}_k
\]
which can be schematized as shown in  \textbf{Algorithm 2}.
\begin{algorithm}[h!]\label{five_GAUSSIAN}
    \caption{ \((I+{\bf V}{\bf \Psi}{\bf V}) {\rho}_k\) Algorithm}%
        \begin{algorithmic}[1]
            \STATE $z_1={\bf V} {\rho}_k$;\\[0.5 mm]
            \STATE $z_2= {\bf \Psi} z_1$;\\[0.5 mm]
            \STATE $z_3={\bf V} z_2$;\\[0.5 mm]
            \STATE $\mathbf q_k=  {\rho}_k + z_3$;
             \end{algorithmic}
\label{PCG}
\end{algorithm}

\noindent The steps 1. and 3. in  \textbf{Algorithm 2} consist in a
matrix-vector product. These products, as detailed in next section,
can be considered discrete Gaussian convolutions and the matrix
${\bf V}$, for one-dimensional state vectors, has Gaussian
structure. Even for state vectors defined on two (or more)
dimensions, the matrix ${\bf V}$ can be represented as product of
two (or more) Gaussian matrices. Since a single matrix-vector
product of this form becomes prohibitively expensive if carried out
explicitly, a computational advantage is gained by employing
Gaussian RFs to mimic the required Gaussian convolution operators.\\
In the previous OceanVar scheme, it was implemented a 1st-RF
algorithm, as described in
\cite{Purser95, Purser}. Here, we study the 3rd-RF introduction, based on
\cite{Young, Vliet}.
\\
The aim of the following sections is to precisely reveal how
the $n$-th order recursive filters are defined and, through the
error analysis, to investigate on their effect in terms of error
estimate and perfomences.

\section{Gaussian recursive filters}
In this section we describe Gaussian recursive filters as
approximations of the discrete Gaussian convolution used in steps 1.
and 3. of \textbf{Algorithm 2}. Let denote by $$g(x)=
\frac{1}{\sigma\sqrt{2\pi}} \exp\left(- \frac{x^2}{2 \sigma
^2}\right)$$ the normalized Gaussian function and by $\bf V$ the
square matrix whose entries are given by
\begin{equation}\label{convGauss5}
{\bf V}_{i,j}=
g(i-j)=\frac{1}{\sigma\sqrt{2\pi}} \exp\left(- \frac{(i-j)^2}{2
\sigma ^2}\right).
\end{equation}
Now let be $s^0=(s_1^0,\ldots,s_m^0)^T$ a vector; the discrete
Gaussian convolution of $s^0$ is a new vector $s=(s_1,\ldots,s_m)^T$
defined by means of the matrix-vector product
\begin{equation}\label{Gauss_oper}
s  =  {\bf V}\otimes s^0\equiv  {\bf V}\;s^0.
\end{equation}
The discrete Gaussian convolution can be considered as a discrete
representation of the continuous Gaussian convolution.  As is well
known, the continuous Gaussian convolution of a function $s^0$ with
the normalized Gaussian function $g$  is a new function $s$ defined
as follows:
\begin{equation}\label{convGauss}
s(x)=[g  \otimes s^0](x) = \int_{-\infty}^{+\infty} g(x- \tau)
s^0(\tau) d \tau.
\end{equation}
Discrete and continuous Gaussian convolutions are strictly related.
This fact could be seen as follows. Let assume that
$$I=\{x_1<x_2<\ldots <x_{m+1}\}$$ is a grid of evaluation points and
let set for $i=1,\ldots,m$
$$s_i\equiv s(x_i),\qquad  s_i^0\equiv s^0(x_i)\quad \text{and}\quad \Delta x_i=
x_{i+1}-x_i=1.$$ 
By  assuming that $s^0$ is $0$ outside of $[x_1,x_{m+1}]$ and by
discretizing the integral (\ref{convGauss}) with a rectangular
 rule, we obtain
\begin{equation*}\label{convGauss2}
s_i=\int_{-\infty}^{+\infty}\!\!\!\!\!\!\!g(x_i- \tau) s^0(\tau) d
\tau=\int_{x_1}^{x_{m+1}}\!\!\!\!\!\!\!g(x_i- \tau) s^0(\tau) d
\tau=\qquad
\end{equation*}
\begin{equation*}\label{convGauss3}
\quad \ =\sum_{j=1}^m \int_{x_j}^{x_{j+1}} \!\!\!\!\!\!\!g(x_i-
\tau) s^0(\tau) d \tau \approx \sum_{j=1}^m  \Delta x_j g(x_i- x_j)
s^0_j=\qquad
\end{equation*}
\begin{equation}\label{convGauss3BIS}
\quad \ = \sum_{j=1}^m  g(i-j) s^0_j = \sum_{j=1}^m  {\bf V}_{i,j}
s^0_j= ({\bf V} s^0)_i.
\end{equation}
An optimal way for approximating the values $s_i$ is given by
Gaussian recursive filters. The $n$-order RF filter computes the
vector $s^K=(s^K_1,\ldots,s^K_m)^T$ as follows:
\begin{equation}\label{forward&backward_n}
\left\{\begin{array}{l}
p_i^k=\beta_i s_i^{k-1}+ \displaystyle{\sum_{j=1}^n} \alpha_{i,j} p_{i-j}^k  \,\,\,\ i=1,\ldots,m\\[0.4cm]
s_i^k=\beta_i p_i^{k}+ \displaystyle{\sum_{j=1}^n} \alpha_{i,j}
s_{i+j}^k  \,\,\,\ i=m,\ldots,1
\end{array}
\right..
\end{equation}
The iteration counter $k$ goes from $1$ to $K$, where $K$ is the
total number of filter iterations. Observe that values $p_1^k,\ldots
p_n^k$ are computed taking in the sums terms $\alpha_{i,j}
p_{i-j}^k$  provided that $i-j\geq 1$. Analogously values
$s_m^k,\ldots s_{m-n+1}^k$ are computed taking  in the sums terms
$\alpha_{i,j} s_{i+j}^k$ provided that $i+j\leq m$. The values \(
\alpha_{i,j}\) and $\beta_{i}$, at each grid point $x_i$,  are often
called \emph{smoothing coefficients} and they obey to the constraint
$$\beta_i=1-\sum_{j=1}^n \alpha_{i,j}.$$ In this paper we deal with
first-order and third-order RFs. The first-order RF expression
($n=1$) becomes:
\begin{equation}\label{forward&backword1D}
\left\{\begin{array}{l} p_1^k=\beta_1 s_1^{k-1}, \\
p_i^k=\beta_i s_i^{k-1}+ \alpha_i p_{i-1}^k  \quad i=2,\ldots,m\\[0.2cm]
s_m^k=\beta_m p_m^{k},  \\ \, s_i^k=\beta_i p_i^{k}+ \alpha_i
s_{i+1}^k \qquad i=m-1,\ldots,1.
\end{array}
\right.
\end{equation}
If $R_i$ is the correlation radius at $x_i$, by setting $$\sigma_i =
\frac{R_i}{\Delta x_i} \quad \text{and} \quad E_i\!=\!\frac{K \Delta
x_i^2}{R_i^2}\!=\!\frac{K}{\sigma_i^2},$$ coefficients $\alpha_i$ e
$\beta_i$ are given by \cite{Purser95}:
\begin{equation}\label{alphai-betai}
\alpha_i\!=\!1+E_i - \sqrt{ E_i(E_i\!+\!2)}, \ \ \ \  \ \
\beta_i
\!=\!\sqrt{ E_i(E_i\!+\!2)}-E_i. \ \
\end{equation}
The third-order RF expression ($n=3$) becomes:
\begin{equation}\label{forward&backward_3}
\left\{\begin{array}{l}
p_i^k=\beta_i s_i^{k-1}+ \displaystyle{\sum_{j=1}^3} \alpha_{i,j} p_{i-j}^k  \,\,\,\ i=1,...,m\\[0.4cm]
s_i^k=\beta_i p_i^{k}+ \displaystyle{\sum_{j=1}^3} \alpha_{i,j}
s_{i+j}^k  \,\,\,\ i=m,\ldots,1.
\end{array}
\right.
\end{equation}
Third-order RF coefficients $\alpha_{i,1},\alpha_{i,2},\alpha_{i,3}$
and $\beta_i$, for one only filter iteration ($K=1$), are computed
in \cite{farina2}. If
\begin{equation*}\label{ai}
a_{i}=3.738128 + 5.788982 \sigma_i + 3.382473 \sigma_i^2+
\sigma_i^3.
\end{equation*}
the coefficients expressions are:
\[ \begin{array}{l}
\alpha_{i,1}= (5.788982 \sigma_i+  6.764946  \sigma_i^2 + 3
\sigma_i^3)/a_i
\\[.2cm]
\alpha_{i,2}=-(3.382473 \sigma_i^2 +3 \sigma_i^3)/a_i\\[.2cm]
\alpha_{i,3}= \sigma_i^3 /a_i\\[.2cm]
\beta_i=1-(\alpha_{i,1}+\alpha_{i,2}+\alpha_{i,3})=3.738128/a_i.
\end{array}
\]
In   \cite{Vliet} is proposed the use
of a value $q=q(\sigma_i)$ instead of $\sigma_i$. The $q$ value is:
\begin{equation}\label{qsigma}
q(\sigma_i)=\!\!\left\{\!\!\begin{array}{l}
          0.98711  \sigma_i - 0.96330 \qquad \text{if} \quad \sigma_i>2.5\\[0.2cm]
          3.97156 - 4.14554 \sqrt{1 - 0.26891\sigma_i} \quad
          \text{oth.}
\end{array}
\right.
\end{equation}
In order to understand how Gaussian RFs approximate the discrete
Gaussian convolution it is useful to represent them in terms of
matrix formulation. As explained in \cite{Dahlquist}, the $n$-order
recursive filter computes $s^K$ from $s^{0}$ as the solution of the
linear system
\begin{equation}\label{sysL-U}
(LU)^Ks^K=s^0,
\end{equation}
where matrices $L$ and $U$ are
respectively lower and upper band triangular with nonzero entries
\begin{equation}\label{L-U}
U_{i,i}=L_{i,i}=\frac{1}{ \beta_i},\qquad
L_{i,i-j}=U_{i,i+j}=-\frac{\alpha_{i,j}}{ \beta_i}.
\end{equation}
By formally inverting the linear system (\ref{sysL-U}) it results
\begin{equation}\label{Filter_oper}
s^K={\bf F_n^{(K)}}\;s^0,
\end{equation}
where ${\bf F_n^{(K)}}\equiv (LU)^{-K}$. A direct expression of
${\bf F_n^{(K)}}$ and its norm could be obtained, for instance, for
the first order recursive filter in the homogenus case
($\sigma_i=\sigma$). However, in the following, it will be shown
that ${\bf F_n^{(K)}}$ has always bounded norm, i.e.
\begin{equation}\label{LIM}
\|{\bf
F_n^{(K)}}\|_\infty\leq1.
\end{equation}
 Observe that ${\bf F_n^{(K)}}$ is the
matrix operator that substitutes the Gaussian operator ${\bf V}$ in
(\ref{Gauss_oper}), then a measure of how well $s^K$ approximates
$s$ can be derived in terms of the operator distance $$\|\bf V-{\bf
F_n^{(K)}}\|_\infty.$$ Ideally one would expect that $\|\bf V-{\bf
F_n^{(K)}}\|$ goes to $0$ (and $s^K \to s$) as $K$ approaches to
$\infty$,  yet this does not happen due to the presence of edge
effects. In the next sections we will investigate about the numerical
behaviour of the distance $\|\bf V-{\bf F_n^{(K)}}\|$ for some case
study and we will show its effects in the CG algorithm.

\section{RF error analysis}
Here we are interested to analyze the error introduced on the
matrix-vector operation at step 5. of \textbf{Algorithm 1}, when the
Gaussian RF is used instead of the discrete Gaussian convolution. As
previously explained, in terms of matrices, this is equivalent to
change the matrix operator, then \textbf{Algorithm 2} can be
rewritten as shown in \textbf{Algorithm 3}.

\begin{algorithm}[h!]\label{five_RF}
    \caption{ \((I+{\bf F_n^{(K)}}{\bf \Psi}{\bf F_n^{(K)}})  \widetilde{\rho}_k\) Algorithm}%
        \begin{algorithmic}[1]
            \STATE $\widetilde{z}_1={\bf F_n^{(K)}}{\widetilde{\rho}}_k$;\\[0.5 mm]
            \STATE $\widetilde{z}_2= {\bf \Psi} \widetilde{z}_1$;\\[0.5 mm]
            \STATE $\widetilde{z}_3={\bf F_n^{(K)}} \widetilde{z}_2$;\\[0.5 mm]
            \STATE $\mathbf{\widetilde{ q}}_k=\widetilde{\rho}_k + \widetilde{z}_3$;
             \end{algorithmic}
\label{PCG2}
\end{algorithm}

\noindent Now we are able to give the main result of this paper:
indeed the following theorem furnishes an upper bound for the error
$\mathbf{q}_k-\mathbf{\widetilde{ q}}_k$, made at each single
iteration $k$ of the CG (\textbf{Algorithm 1}). This bound involves
the operator norms
$$\Vert{\bf F_n^{(K)}}\Vert_\infty, \quad \Vert{\bf \Psi
}\Vert_\infty, \quad \Vert{\bf V}\Vert_\infty,$$ the distance
$\Vert\bf V-{\bf F_n^{(K)}}\Vert_\infty$ and the error
${{\rho}}_k-{\widetilde{\rho}}_k$ accumulated on ${{\rho}}_k$
at previous iterations.\\
\noindent  \begin{theorem} Let be ${{\rho}}_k$, $\widetilde{\rho}_k$,
$\mathbf{ q}_k$, $\mathbf{\widetilde{ q}}_k$ as in \textbf{Algorithm
2} and \textbf{Algorithm 3}.  Let be $\Vert \cdot \Vert=\Vert \cdot \Vert_\infty $ and
let denote by
$$e_k={{\rho}}_k-\widetilde{{\rho}}_k$$ the difference between
values ${{\rho}}_k$ and  $\widetilde{{\rho}}_k$. Then it holds
\begin{eqnarray}\label{theor-bound}
\Vert\mathbf{{ q}}_k-\mathbf{\widetilde{ q}}_k\Vert
\leq \left(1+ \Vert {\bf V}\Vert\!\cdot\!  \Vert{\bf \Psi}\Vert \!\cdot \! \Vert{\bf
V}\Vert \right)\! \cdot \! \Vert e_k\Vert +  \ \qquad \ \ \  \ \ \notag
\\[.2cm]
\ \ + \Vert {\bf F_n^{(K)}} \!-\!{\bf V} \Vert \!\cdot\!
\Vert{\bf \Psi}\Vert  \!\cdot \!\big(\Vert{\bf V} \Vert \!+\!\Vert {\bf F_n^{(K)}}\Vert\big)
\!\cdot\! \Vert\widetilde{\rho}_k\Vert.
\end{eqnarray}
\textbf{Proof:} A direct proof follows by using the values $z_i$ and $\widetilde{z}_i$ introduced in
 Algorithm 2 and in  Algorithm 3. It holds:
\begin{eqnarray}\label{proof-1}
\Vert{z}_1\!-\!\widetilde{z_1}\Vert\!=\!\Vert{{\bf V}\rho _k \!-\!\bf F_n^{(K)}}\widetilde{\rho }_k
\Vert\!=\!
\Vert{\bf V}\rho _k \!-\! {\bf V}\widetilde{\rho}_k\!+\!{\bf V}\widetilde{\rho}_k\!-\!{\bf F_n^{(K)}}
\widetilde{\rho}_k
\Vert\!\leq\notag \\
\leq
\Vert{\bf V} {\rho}_k\!-\!
{\bf V}\widetilde{\rho}_k\Vert
\!+\!\Vert{{\bf V} \widetilde{\rho}_k\!-\!\bf F_n^{(K)}}\widetilde{\rho }_k\Vert\!\leq\!
\notag \hspace{2.46cm}\\
\leq \Vert{\bf V}\Vert  \cdot  \Vert e_k\Vert+\Vert{\bf V}-{\bf F_n^{(K)}} \Vert
\cdot  \Vert \widetilde{\rho }_k \Vert. \notag \hspace{2.63cm}
\end{eqnarray}
Then, for the difference ${z}_2-\widetilde{z}_2$, we get the bound
\begin{eqnarray}\label{proof-2}
\Vert{z}_2-\widetilde{z}_2\Vert = \Vert{\bf \Psi} {z_1}- {\bf \Psi} \widetilde{z_1} \Vert\leq \Vert{\bf \Psi} \Vert
\cdot \Vert{z}_1\!-\!\widetilde{z}_1\Vert
\leq \notag \hspace{1.7cm}\\
\leq \Vert{\bf \Psi} \Vert \cdot   \Vert{\bf V}\Vert  \cdot  \Vert e_k\Vert + \Vert{\bf \Psi} \Vert   \cdot
\Vert{\bf V}-{\bf F_n^{(K)}}\Vert   \cdot   \Vert \widetilde{\rho }_k \Vert. \hspace{0.15cm}\notag
\end{eqnarray}
Hence, for the difference ${z}_3-\widetilde{z}_3$, we obtain
\begin{eqnarray}\label{proof-3}
\Vert{z}_3\!-\!\widetilde{z}_3\Vert \!=\!
\Vert{\bf V} {z_2}\!-\! {\bf F_n^{(K)}}\widetilde{z_2}\Vert\!=\!
\Vert{\bf V} {z_2}\!-\! {\bf V}\widetilde{z}_2\!+\! {\bf V}\widetilde{z}_2\!-\!
 {\bf F_n^{(K)}}\widetilde{z_2}\Vert\!\leq \notag \hspace{.7cm} \\ 
\leq \Vert{\bf V} \Vert \cdot \Vert {z_2}- \widetilde{z_2}\Vert+\Vert{\bf V}-{\bf F_n^{(K)}} \Vert
\cdot \Vert\widetilde{z_2}\Vert \leq \notag \hspace{2.3cm}\\
\leq \Vert{\bf V} \Vert \!\cdot\! \Vert {z_2}\!-\! \widetilde{z_2}\Vert\!+\!\Vert{\bf V}\!-\!{\bf F_n^{(K)}}
\Vert
\!\cdot \!\Vert{\bf \Psi} \Vert\! \cdot\! \Vert{\bf F_n^{(K)}} \Vert \cdot \Vert
\widetilde{\rho }_k \Vert \leq \notag \hspace{.7cm} \\
\leq \Vert{\bf V}\Vert  \!\cdot  \!\Vert{\bf \Psi} \Vert \!\cdot   \!
\Vert{\bf V}\Vert  \!\cdot  \!\Vert e_k\Vert \!+\!
\Vert{\bf V}\Vert \!\cdot  \!\Vert{\bf \Psi} \Vert   \!\cdot\!
\Vert{\bf V}\!-\!{\bf F_n^{(K)}}\Vert   \!\cdot   \! \Vert \widetilde{\rho }_k \Vert \notag \hspace{.75cm}  \\
 + \Vert{\bf V}\!-\!{\bf F_n^{(K)}}
\Vert
\!\cdot \!\Vert{\bf \Psi} \Vert\! \cdot\! \Vert{\bf F_n^{(K)}} \Vert \cdot \Vert
\widetilde{\rho }_k \Vert =\notag \hspace{3.3cm} \\
\Vert{\bf V}\Vert  \!\cdot  \!\Vert{\bf \Psi} \Vert \!\cdot   \!
\Vert{\bf V}\Vert  \!\cdot  \!\Vert e_k\Vert \!+\!
\Vert{\bf V}\!-\!{\bf F_n^{(K)}}\Vert   \!\cdot
\!\Vert{\bf \Psi} \Vert   \big(\Vert{\bf V} \Vert \!+\!\Vert {\bf F_n^{(K)}}\Vert\!\big)
\!\cdot \!\Vert
\widetilde{\rho }_k \Vert \notag    \hspace{.7cm}
\end{eqnarray}
In the second-last inequality we used the fact that
$$\Vert\widetilde{z}_2\Vert=
\Vert{\bf \Psi} \widetilde{z}_1\Vert=
\Vert{\bf \Psi} {\bf F_n^{(K)}}\widetilde{\rho }_k \Vert\leq
\Vert{\bf \Psi} \Vert \cdot \Vert{\bf F_n^{(K)}} \Vert \cdot \Vert
\widetilde{\rho }_k \Vert.$$
Finally, observing that
\begin{eqnarray}\label{proof-4}
\Vert\mathbf{{ q}}_k - \mathbf{\widetilde{ q}}_k\Vert
 = \Vert {\rho}_k  +   {z}_3 - (\widetilde{\rho}_k  +  \widetilde{z}_3)\Vert \leq  \notag \hspace{2cm}\\
\leq \Vert {\rho}_k  - \widetilde{\rho}_k\Vert +  \Vert z_3 - \widetilde{z}_3\Vert
 =
\Vert e_k \Vert +\Vert z_3-\widetilde{z}_3\Vert, \notag\hspace{-.7cm}
\end{eqnarray}
and taking the upper bound of $\Vert z_3-\widetilde{z}_3\Vert$, the thesis is proved.\\
\phantom{c}\hfill $\diamond$
\end{theorem}

Previous theorem shows that, at each iteration of the CG algorithm,
the error bound on  the computed value \({\bf {q}_k}\)  at step 5., is characterized by two main terms: 
the first term can be considered as 
the contribution of the standard forward error analysis and it is not significant, if $\|e_k\|$ is small; the second term
highlights  the effect of the introduction  of the RF.
More in detail, at each iteration step, the computed value \({\bf {q}_k}\) is biased by a quantity 
proportional to three factors:
\begin{itemize}
\item the distance between the original operator 
(the Gaussian operator {\bf V}) and its approximation (the operator \({\bf F^{(K)}_n}\));
\item the norm of \({\bf \Psi}\);
\item the sum of the operator norms \(\|{\bf F^{(K)}_n}\|\) and \(\|{\bf V}\|\).
\end{itemize}

\begin{center}
\textbf{Table 1}: Operator norms \\[0.2cm]
\begin{tabular}{ ||l || c || c || }\hline\hline
  $\sigma$ & $||\bf F^{(1)}_1 ||_\infty$ & $||\bf F^{(1)}_3||_\infty$ \\ \hline
 5  &  0.9920  &  0.9897\\ \hline
 20 &  0.9012  &  0.8537\\ \hline
 50 &  0.9489  &  0.8950\\ \hline
 \end{tabular}
 \end{center}

As shown in (\ref{LIM}) the norm \({\bf \Psi}\) is bounded. Besides, the norm of  {\bf V} is always less or equal to one 
(because it comes from the discretization of the of the continuous Gaussian convolution). The 
norm of  \({\bf F^{(K)}_n}\) is bounded by one too. This fact can be seen by observing the Table 1, where
 we consider several tests by varying data distributions in the homogeneous case ($\sigma_i=\sigma$),
for 1st-RF and 3rd-RF. Starting from these considerations, the error estimate of \emph{Theorem 4.1} can be 
spcialized as:
\begin{eqnarray}\label{theor-bound_final}
\Vert\mathbf{\widetilde{ q}}_k-\mathbf{{ q}}_k\Vert
\leq \left(1\!+\! \Vert{\bf \Psi}\Vert \right) \Vert e_k\Vert
\!+ 2 \Vert {\bf F_n^{(K)}} \!-\!{\bf V} \Vert \cdot
\Vert{\bf \Psi}\Vert
\Vert\rho_k\Vert.
\end{eqnarray}

\section{Experimental Results}
In this section we report some experiments to confirm the discussed theoretical results. In the first part, we deal
with the approximations of the discrete operator $\bf V$  with the first order and  of the third order 
$\bf F^{(K)}_1$ and $\bf F^{(1)}_3$ respectively. In the last subsection, we analyze the improving 
in the performance and in the accuracy terms of the third order RF applied to the case study.

\subsection{1st-RF and 3rd-RF operators}
\begin{figure}[ht!]
\centering
        \includegraphics[width=0.30\textwidth]{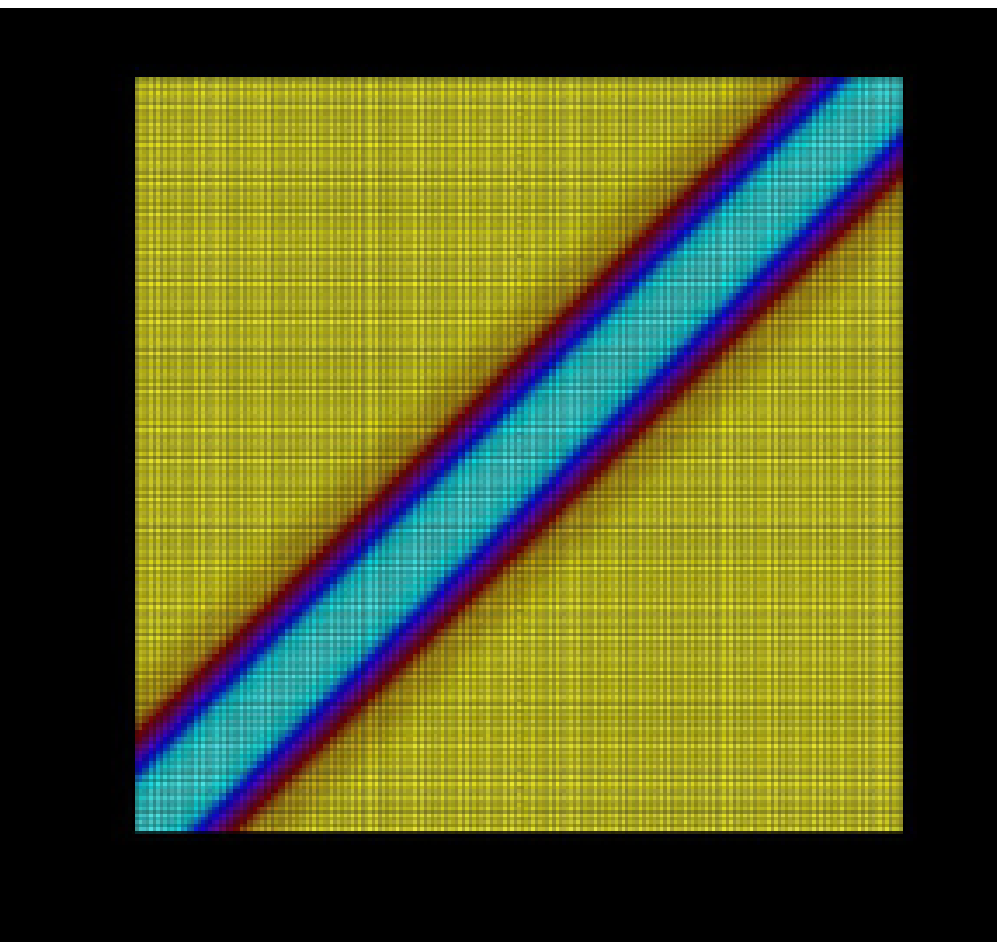}
         \includegraphics[width=0.30\textwidth]{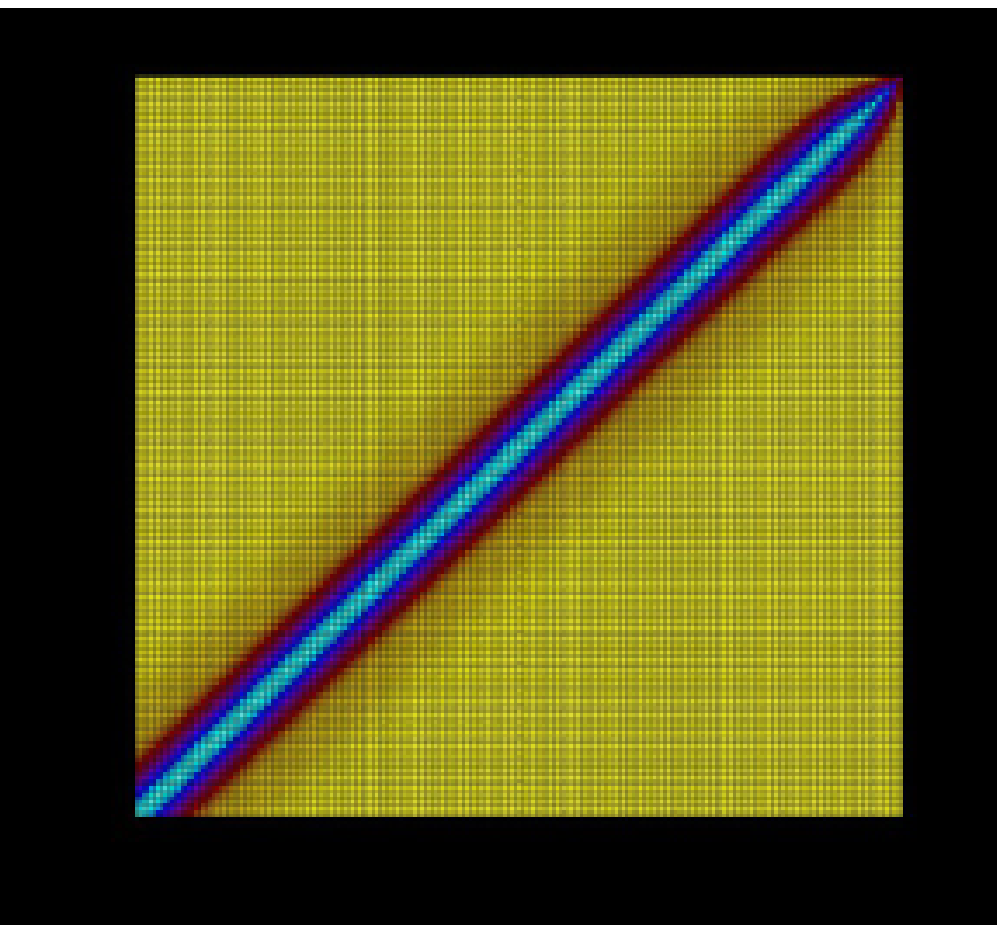}
\caption{
\small{ {\bf Top}. Discrete Gaussian convolution operator $\bf V$. {\bf Bottom}. 1-st order
recursive filter operator $F_1$   }
}
\label{fig:OpG1}
\end{figure}

\begin{figure}[ht!]
\centering
          \includegraphics[width=0.30\textwidth]{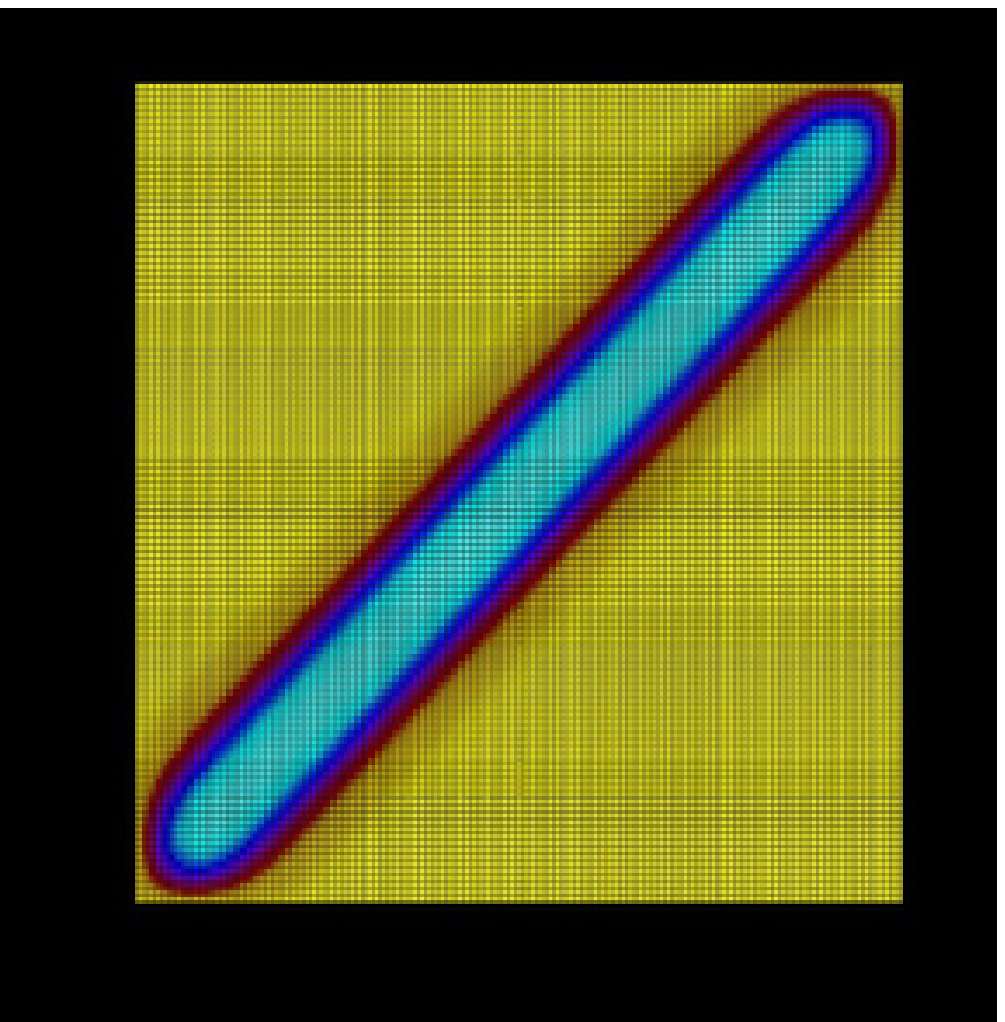}
           \includegraphics[width=0.30\textwidth]{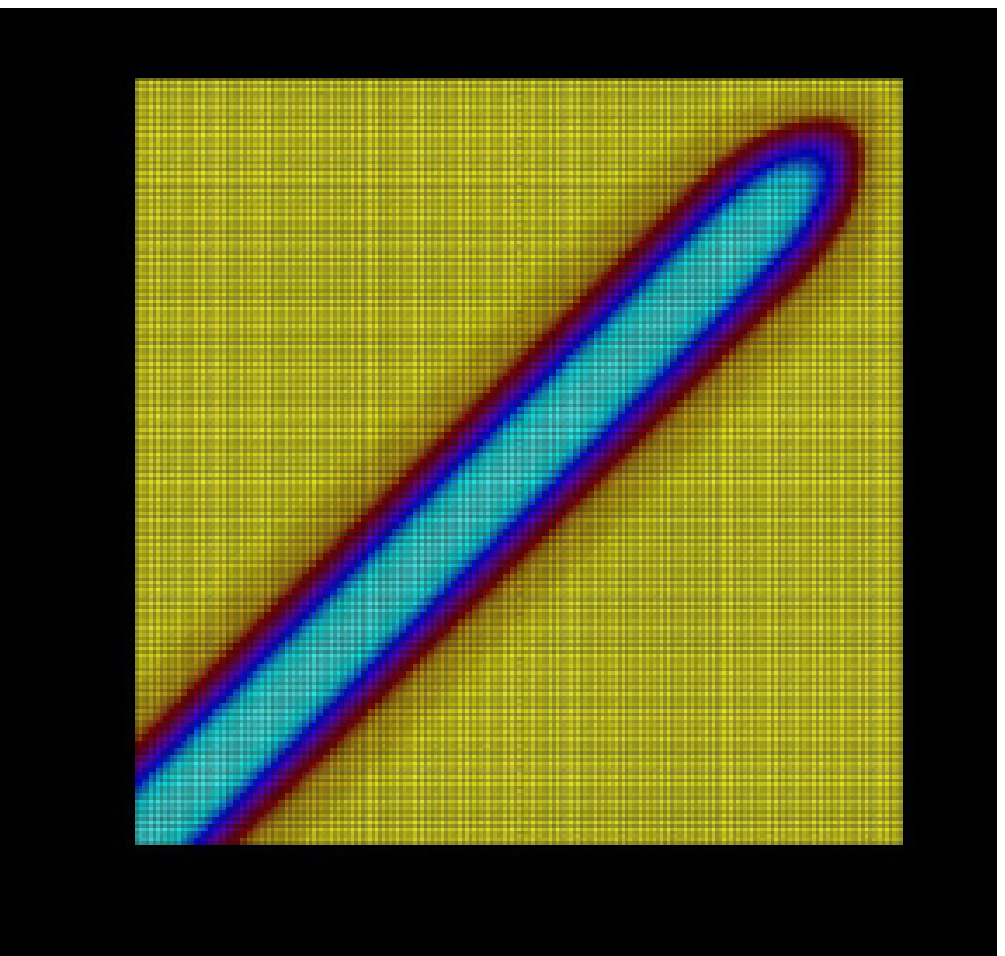}
\caption{
\small{ {\bf Top}. 1-st order recursive filter operator $F^{(50)}_1$  with $50$ iterations. {\bf Bottom}. 3-rd order
recursive filter operator $F^{(1)}_3$   }
}
\label{fig:OpF3}
\end{figure}

In the following experiments, we construct the operators $\bf V$, $\bf F^{(1)}_1$, $\bf F_1^{(50)}$ and $\bf F^{(1)}_3$  
in the  case of $m=601$ samples of a random vector $\bf s^0$. We assume that $\bf s^0$ comes 
from a uniform grid with   homogeneous condition $\sigma_i=\sigma=15$.
In Figure \ref{fig:OpG1},  it is highlighted that the involved discrete operators have different structures.
In particular, a first qualitative remark is that the operator $\bf F^{(1)}_1$ is a poor approximation of $\bf V$.
Conversely, the operator $\bf F_1^{(50)}$ (Figure \ref{fig:OpF3} on the top)  is very close to $\bf V$ but, as for $\bf F^{(1)}_1$,  there 
are significant differences with $\bf V$ in the bottom left and in the top right corners. These dissimilarities in the edges, by a numerical point of view,  
give some kind of artifacts in the computed convolutions, that determine a vector $\bf s$ with  components, in the initial and final positions,
that decay to zero.

Figure \ref{fig:OpF3} bottom shows that the operator $\bf F^{(1)}_3$  is closer then $\bf F^{(1)}_1$ and $\bf F^{(50)}_1$ to the discrete convolution $\bf V$. In particular,  this recursive filter is able to reproduce $\bf V$  more accurately in the bottom left corner, but unfortunately it  does not give  good results  on top right corner. In Table 2, for random distributions with homogeneous condition ($\sigma_i=\sigma$),
we underline the edge effects by measuring the norms between the discrete convolution $\bf V$ and the RF filters. Although the
$||\bf F_n^{(K)} - V||_\infty$ ideally goes to zero as  $k$ goes to $+\infty$, this does not happen in practice as observed below.
\begin{center}
\textbf{Table 2}: Distance metrics \\[0.2cm]
\begin{tabular}{ ||l || c || c || c || }\hline\hline
$\sigma$ & $||\bf F^{(1)}_1 - V||_\infty$ & $||\bf F^{(50)}_1 - V||_\infty$ & $||\bf F^{(1)}_3 - V||_\infty$ \\ \hline
 5  &  0.2977  &  0.3800  &  0.5346 \\ \hline
 10  &  0.3895  &  0.4397  &  0.5890 \\ \hline
 25  & 0.4533   & 0.4758   & 0.6221 \\ \hline
 50  & 0.4686  &  0.4809   & 0.6125 \\ \hline
  \end{tabular}
  \end{center}

\begin{figure*}[ht!]
\centering
          \includegraphics[width=1\textwidth]{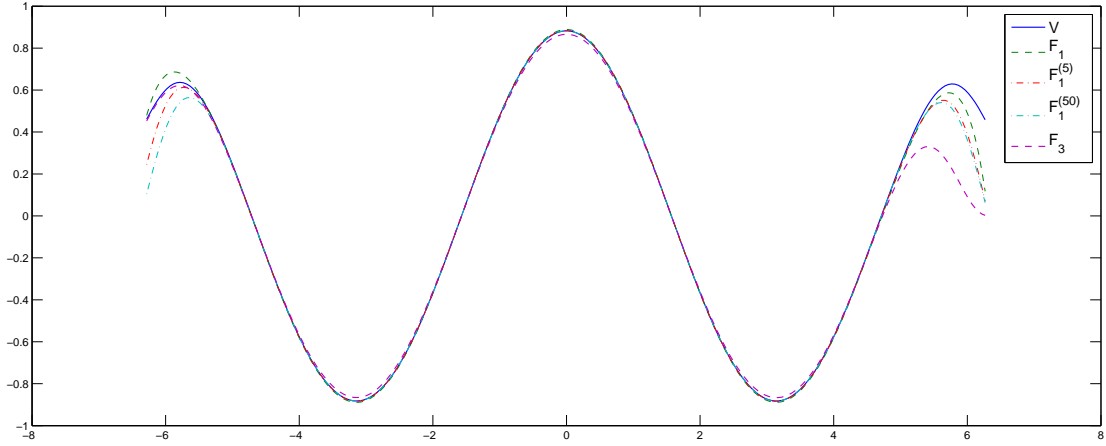}
\caption{
\small{Discrete convolution $\bf V$ and Gaussian recursive filtering $\bf F^{(K)}_1$ with $1,5,50$ iterations and $\bf F^{(1)}_3$ applied to
$n=252$ samples of the periodic function $s^0=\cos (x)$ in $[-2\pi, 2\pi]$. }
}
\label{fig:OpF4}
\end{figure*}

In order to bring out these considerations, we show the application of $\bf V$, $\bf F^{(K)}_1$ and $\bf F^{(1)}_3$  to a periodic signal $\bf s^0$. We choose $m=252$ samples of the $\cos$ function in $[-2\pi, 2\pi]$ and we perform simulations by using the 1-st RF with  $1, 5$ and $50$ iterations
and 3-rd RF with  one iteration.
In Figure \ref{fig:OpF4} it is shown the computed Gaussian convolution and the poor approximation of ${\bf V}s^0$ on the right side of the test interval, due to
the edge effects. A nice result is that our $\bf F^{(1)}_3$ convolution operator gives  better results on the left side of the domain. \\[0.2cm]
Finally, we give some considerations about the accuracy of the studied Gaussian RF schemes,  when they are applied to
the Dirac rectangular impulse\\[0.2cm]
\centerline{$s^0=(0,\ldots,0,1,0,\ldots).$}\\[0.2cm]
We choose a one-dimensional  grid of \(m=301\) points, a constant correlation radius \(R=120,km\) , a constant grid
space \( \Delta x =6\,km\) and  $\sigma=R/ \Delta x=20$.
In the numerical experiments to avoid the edge effects, we only consider  $\bar{m}=221$ central values of $s^K$, i.e.\\[0.2cm]
\centerline{$\bar{s}^K=(s^K_{2\sigma},s^K_{2\sigma+1},\ldots,s^K_{m-2\sigma-1},s^K_{m-2\sigma})$.}\\[0.2cm]
Similarly, in  Table 3  we measure the operator distances we use $\bf ||\bar{F}^{(1)}_1-\bar{V}||_\infty$  and  $\bf ||\bar{F}^{(1)}_3-\bar{V}||_\infty$,
where $\bar{V}$, $\bar{F}_1^{(1)}$ and $\bar{F}_3^{(1)}$ indicate the submatrices obtained, neglecting  first and last $2\sigma-1$ rows and columns.\\
\newpage
\begin{center}
\textbf{Table 3}: Convergence history\\[0.2cm]
\begin{tabular}{ ||l || c || c || }\hline\hline
K & $\bf ||\bar{F}^{(K)}_1-\bar{V}||_\infty$ & $\bf ||\bar{F}^{(K)}_3-\bar{V}||_\infty$ \\ \hline
 1 & 0.211 & \textbf{0.0424} \\ \hline
  2 & 0.13 &  --\\ \hline
  5 & 0.078 & -- \\ \hline
  50 & 0.048 & --\\ \hline
  100 & 0.0429 & -- \\ \hline
  500 & 0.0414 & -- \\ \hline
  \end{tabular}
  \end{center}
 
\begin{figure}[h!]
\centering
       \includegraphics[width=0.85\hsize]{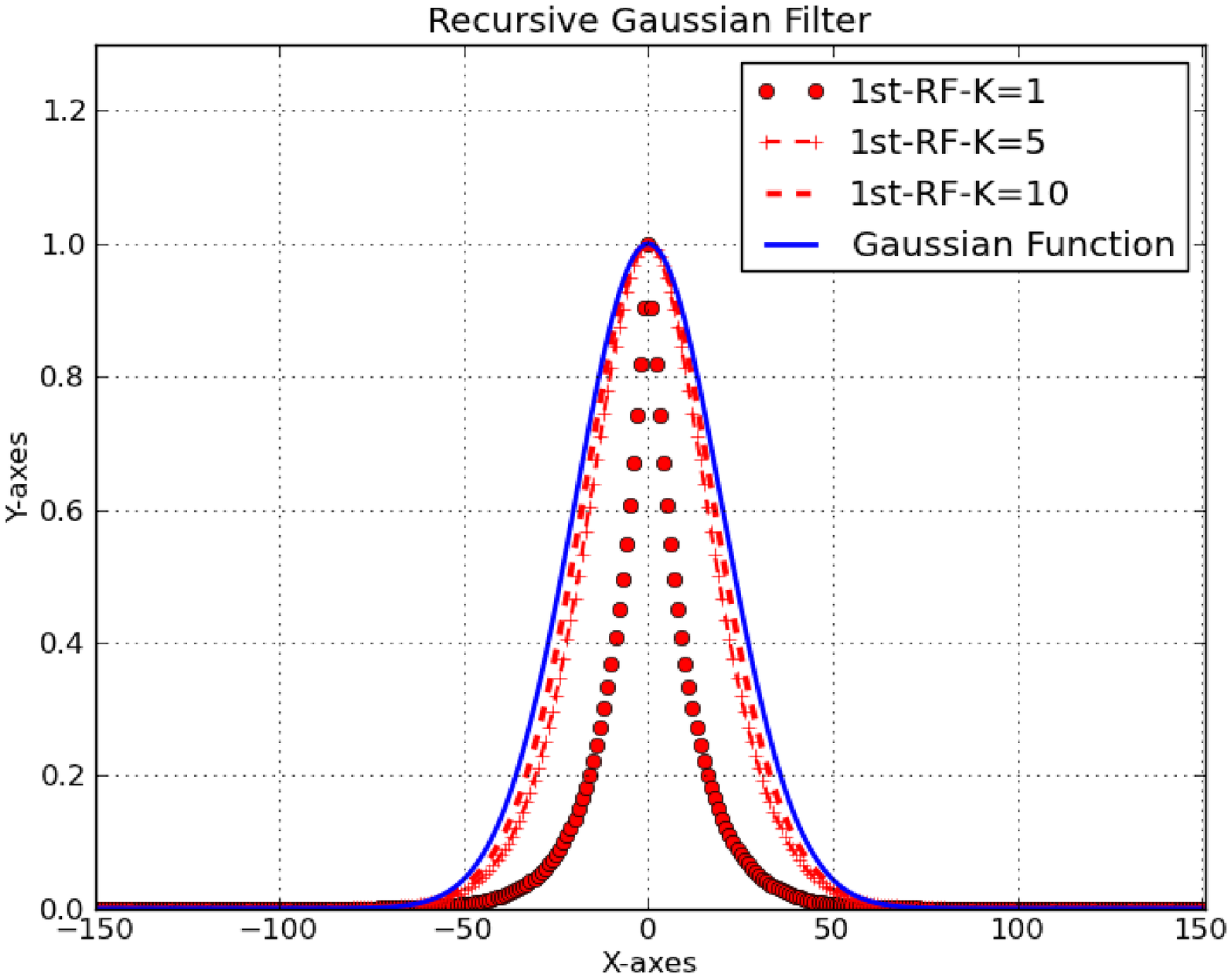}
        \includegraphics[width=0.85\linewidth]{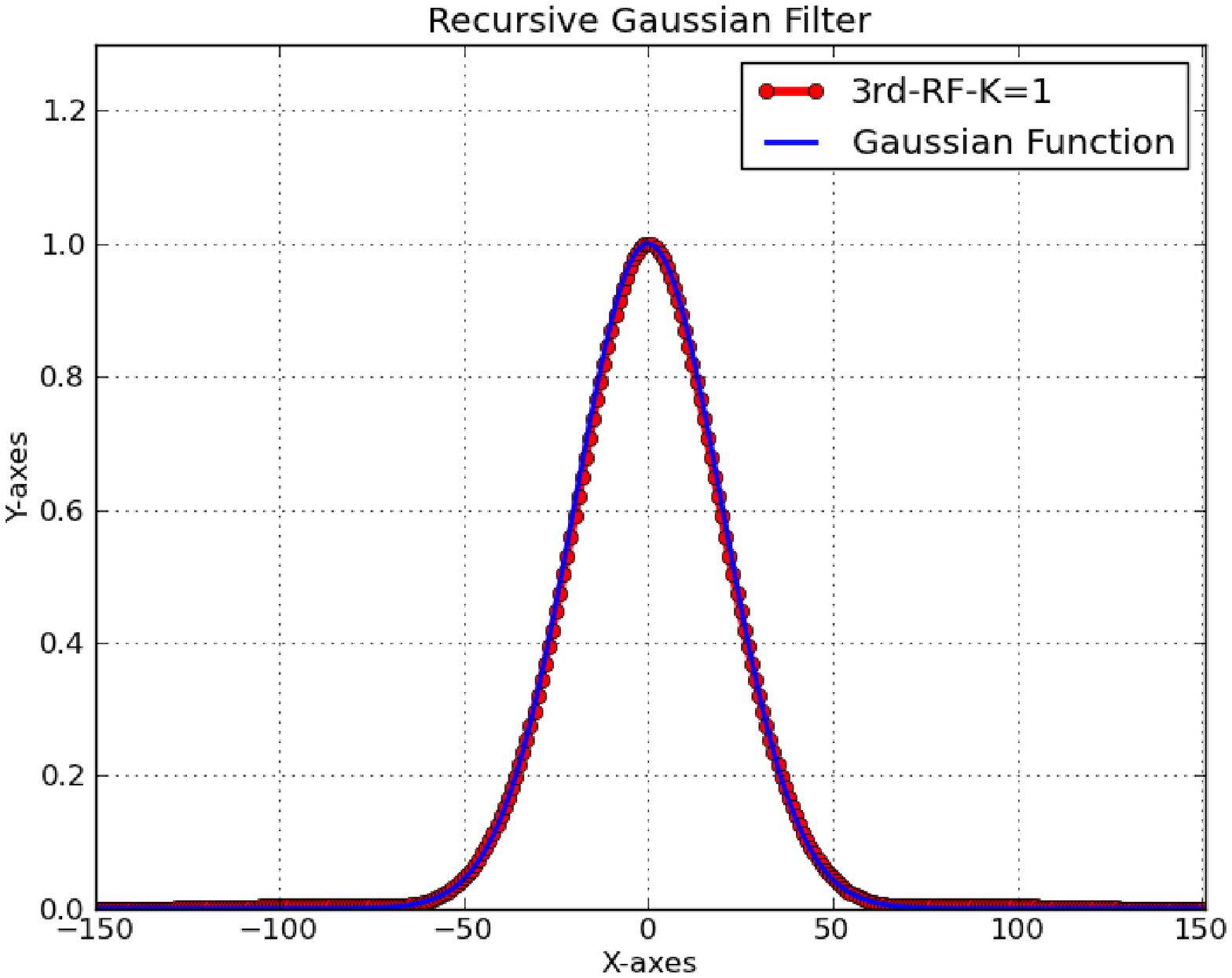}

\caption{
\small{{\bf Top}.   The discrete Gaussian convolution ${\bf {V}}s^0$ (blue) and ${\bf {F}^{(K)}_1} s^0$ for  $K=1,5,10$ (red). {\bf Bottom} The discrete Gaussian convolution ${\bf {V}}s^0$ (blue) and ${\bf {F}^{(1)}_3} s^0$  (red).}
}
\label{fig:pr}
\end{figure}

These case study shows that , neglecting the edge effects, the 3-rd RF filter is more accurate the the 1st-RF order with few iterations. This fact is 
evident by observing  the results in Figure \ref{fig:pr} and the operator norms in Table 3. Finally, we remark that the 1-st order RF has 
to use 100 iteration in order to obtain the same accuracy of the 3-rd order RF. This is a very interesting numerical feature of the third order filter.

\subsection{A case study: Ocean Var}
The theoretical considerations of the previous sections are useful to understand the accuracy improvement in the
real experiments on Ocean Var. The preconditioned  CG is a numerical kernel intensively used in the model minimizations.
Implementing a more accurate convolution operators gives benefits on the convergence of GC and on the overall
data assimilation scheme \cite{farina2} .
Here we report experimental results of the 3rd-RF in a  Global Ocean implementation of OceanVar that follows \cite{Storto}.
These results are extensively discussed in the report \cite{farina2}.
In real scenarios \cite{Gall1, farina1}  scientific libraries and  an high performance computing environments are needed.
The case study simulations were carried-out on an IBM cluster using 64 processors. The model resolution was about $1/4$ degree and the horizontal grid was tripolar, as described in \cite{Madec}. This configuration of the model was used at CMCC for global ocean physical reanalyses applications (see \cite{Ferry}). The model has 50 vertical depth levels. The three-dimensional model grid consists of 736141000 grid-points. The comparison between the 1st-RF and 3rd-RF was carried out for a realistic case study, where all in-situ observations of temperature and salinity from Expendable bathythermographs (XBTs), Conductivity, Temperature, Depth (CTDs) Sensors, Argo floats and Tropical mooring arrays were assimilated. The observational profiles are collected, quality-checked and distributed by \cite{Coriolis}.
The global application of the recursive filter accounts for spatially varying and season-dependent correlation length-scales (CLSs). Correlation length-scale were calculated by applying the approximation given in \cite{Belo} to a dataset of monthly anomalies with respect to the monthly climatology, with inter-annual trends removed.



The analysis increments from a 3DVAR applications that uses the 1st-RF with 1, 5 and 10 iterations and the 3rd-RF are shown in Figure 5 with a zoom in the same area of Western Pacific Area as in Figure 5, for the temperature at 100 m of depth. The Figure also displays the differences between the 3rd-RF and the 1st-RF with either 1 or 10 iterations. The patterns of the increments are closely similar, although increments for the case of 1st-RF (K=1) are generally sharper in the case of both short (e.g. off Japan) or long (e.g. off Indonesian region) CLSs. The panels of the differences reveal also that the differences between 3rd-RF and the 1st-RF (K=10) are very small, suggesting once again that the same accuracy of the 3rd-RF can be achieved only with a large number of iterations for the first order recursive filter. Finally, in [ARXIV] was also observed that the 3rd-RF compared to the 1st-RF (K=5) and the 1st-RF (K=10) reduces the wall clock time
of the software respectively of about 27\% and 48\%.

\begin{figure*}[h!]
\centering
        \includegraphics[width=1\hsize]{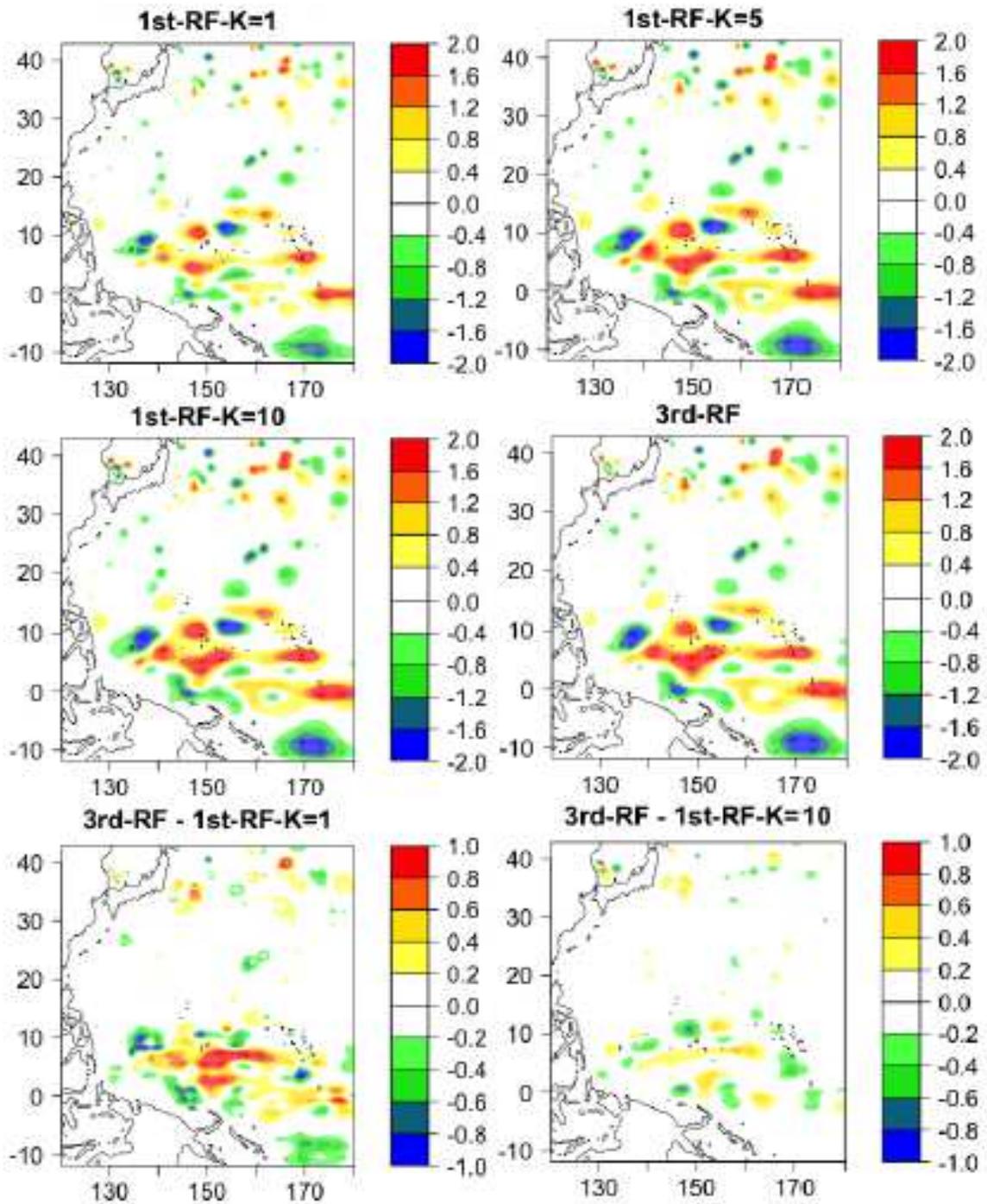}

\caption{
\small{Analysis increments of temperature at 100 m of depth for the Western Pacific for different configurations of the recursive filter (first two rows of panels). Differences of 100 m temperature analysis increments between 3rd-RF and 1st-RF (K=1) and between 3rd-RF and 1st-RF (K=10) (bottom panels).}}
\label{fig:terza}
\end{figure*}

\section{Conclusions}
Recursive Filters (RFs) are a
well known way to approximate the Gaussian convolution and are
intensively applied in the meteorology, in the oceanography and
in forecast models. In this paper, we deal with the oceanographic
3D-Var scheme OceanVar. 
The computational kernel of the OceanVar software is a linear
system solved by means of the Conjugate Gradient (GC)
method. The iteration matrix is related to an error
covariance matrix, with a Gaussian correlation structure.
In other words, at each iteration, a Gaussian convolution is required. 
Generally, this convolution is approximated by a first order RF. 
In this work, we introduced a 3rd-RF filter and we investigated 
about the main sources of error due to
the use of 1st-RF and 3rd-RF operators.
Moreover, we studied how these errors influence the CG algorithm 
and we showed that the third order operator is more
accurate than the first order one.
Finally, theoretical
issues were confirmed by some numerical experiments and
by the reported results in the case study of the OceanVar software.

\end{document}